\documentclass[12pt]{amsart}

\usepackage[english]{babel}

\usepackage{amsmath,amsfonts,amssymb,amsthm}
\usepackage[english]{babel}

\newtheorem{thm}{Theorem}
\newtheorem{cor}[thm]{Corollary}

\theoremstyle{definition}
\newtheorem{defn}{Definition}
\theoremstyle{remark}
\newtheorem{rem}[thm]{Remark}
\numberwithin{equation}{section}
\usepackage{graphicx}

\DeclareGraphicsExtensions{.pdf,.png,.jpg}
\begin{document}

\title{{{TRIANGLE CONICS AND CUBICS}}}

\maketitle

\begin{center}
	Ruslan Skuratovskii\\
	The National Technical University of Ukraine Igor Sikorsky Kiev Polytechnic
	Institute, lecturer \it e-mail: r.skuratovskii@kpi.ua, ruslcomp@mail.ru \\
ORCID: 0000-0002-5692-6123
\end{center}

\begin{center}
	Veronika Strarodub\\
	American University in Bulgaria, student 
\end{center}

\begin{abstract}
 This is a paper about triangle cubics and conics in classical geometry with elements of projective geometry. In recent years, N.J. Wildberger has actively dealt with this topic using an algebraic perspective. Triangle conics were also studied in detail by H.M. Cundy and C.F. Parry recently. The main task of the article was to develop an algorithm for creating curves, which pass through triangle centers. During the research, it was noticed that some different triangle centers in distinct triangles coincide. The simplest example: an incenter in a base triangle is an orthocenter in an excentral triangle. This was the key for creating an algorithm. Indeed, we can match points belonging to one curve (base curve) with other points of another triangle. Therefore, we get a new intersting geometrical object. During the research were derived number of new triangle conics and cubics, were considered their properties in Euclidian space. In addition, was discussed corollaries of the obtained theorems in projective geometry, what proves that all of the descovered results could be transfered to the projeticve plane.
\end{abstract}

\small{Key--Words:  triangle cubics, conics, curves, projective geometry, Euclidian space.}



\section{Introduction}
Centers of triangle and central triangles were studied by  Clark Kimberling \cite{Klimberg}. We cosider curves which pase throw incenter in a base triangle is an orthocenter in an excentral triangle and other triangle curves. Also we studied some properties of Jerabek hyperbola for the mid-arc triangle.

\section{Main result}
 \indent Firstly, we will consider excentral triangle. Correspondence between points of the excentral and base triangles will give us significant results in developing new triangle curves. Below you may observe correspondence table between points of the base and excentral triangles.

\begin{table}[h]
	\begin{center}
    \begin{tabular}{ | l | l |}
	\hline
	\textbf{Base triangle} &  \textbf{Excentral triangle}    \\ \hline
    \textit{I} (incenter) & \textit{H} (orthocenter)\\ \hline
	\textit{O} (circumcenter)  & \textit{E} (nine-point center)  \\ \hline
    \textit{Be} (Bevan point) & \textit{O} (circumcenter) \\ \hline
    \textit{Mi} (mittenpunkt) &  \textit{Sy} (Lemoine point)    \\ \hline
   \parbox[c]{6cm}{ \textit{Mi $'$} (isogonal conjugate of the mittenpunkt   with respect to the base triangle)}   &    \textit{M} (centroid)  \\ \hline
    \textit{Sp} (Speaker point)     &   \textit{Ta} (Taylor point)   \\ \hline
   \textit{Sy} (Lemoine point)    & \parbox[c]{6cm}{ \textit{Sy($H_{1}H_{2}H_{3}$)} (Lemoine point of the orthic triangle)  }     \\ \hline
   \parbox[c]{6cm}{ \textit{Mi$''$} (isogonal conjugate point of the mittenpunkt  with respect to the excetntral triangle)}  & \parbox[c]{6cm}{\textit{GOT} (homothetic center k of the orthic and tangent triangles)}    \\ \hline
    \end{tabular}
    \end{center}
    \indent
    \label{tabular:table1}
    \caption{Correspondence table between points of the base and excentral triangles}
\end{table}

\indent All of the above facts could be easily proved by basic principles of classiical geometry \cite{Klimberg}.
\indent Hence, we may apply derived results for creating new triangle cubics and conics. Firstly Jerabek hyperbola was considered.

\begin{defn}
	{\rm Jerabek hyperbola is a curve which passes through verticies of trianlge, circumcenter, orthocenter, Lemoine point, isogonal conjugate of the de Longchamps point \cite{hyperbola}. }
\end{defn}
     \indent We may observe that Jerabek hyperbola for the extcentral tringle has number of points which corresponde to other ones in the base triangle. The study of such matches gave us significant results\ref{tabular:table2}.

\begin{table}[h]
\begin{center}
\begin{tabular}{ | l | l |}
\hline
\textbf{Jerabek hyperbola for excentral triangle} & \textbf{New hyperbola for the base triangle}\\\hline
 $A, B, C$ (vertices of the base tringle) & $I_{1}$, $I_{2}$, $I_{3}$ (excenters)\\\hline
$O$ (circumcenter) & $Be$ (Bevan point)\\\hline
$H$ (orthocenter)& $I$ (incenter)\\\hline
$Sy$ (Lemoine point) & $Mi$ (mittenpunkt)\\\hline
$L'$ (isogonal conjugate of the de Longchamps point) & $L$ (de Longchamps point)\\\hline
 \end{tabular}
\end{center}
\indent
\label{tabular:table2}
\caption{Mtaching points for Jerabek hyperbola}
\end{table}

\begin{figure}[h]
	\center{\includegraphics[scale=0.5]{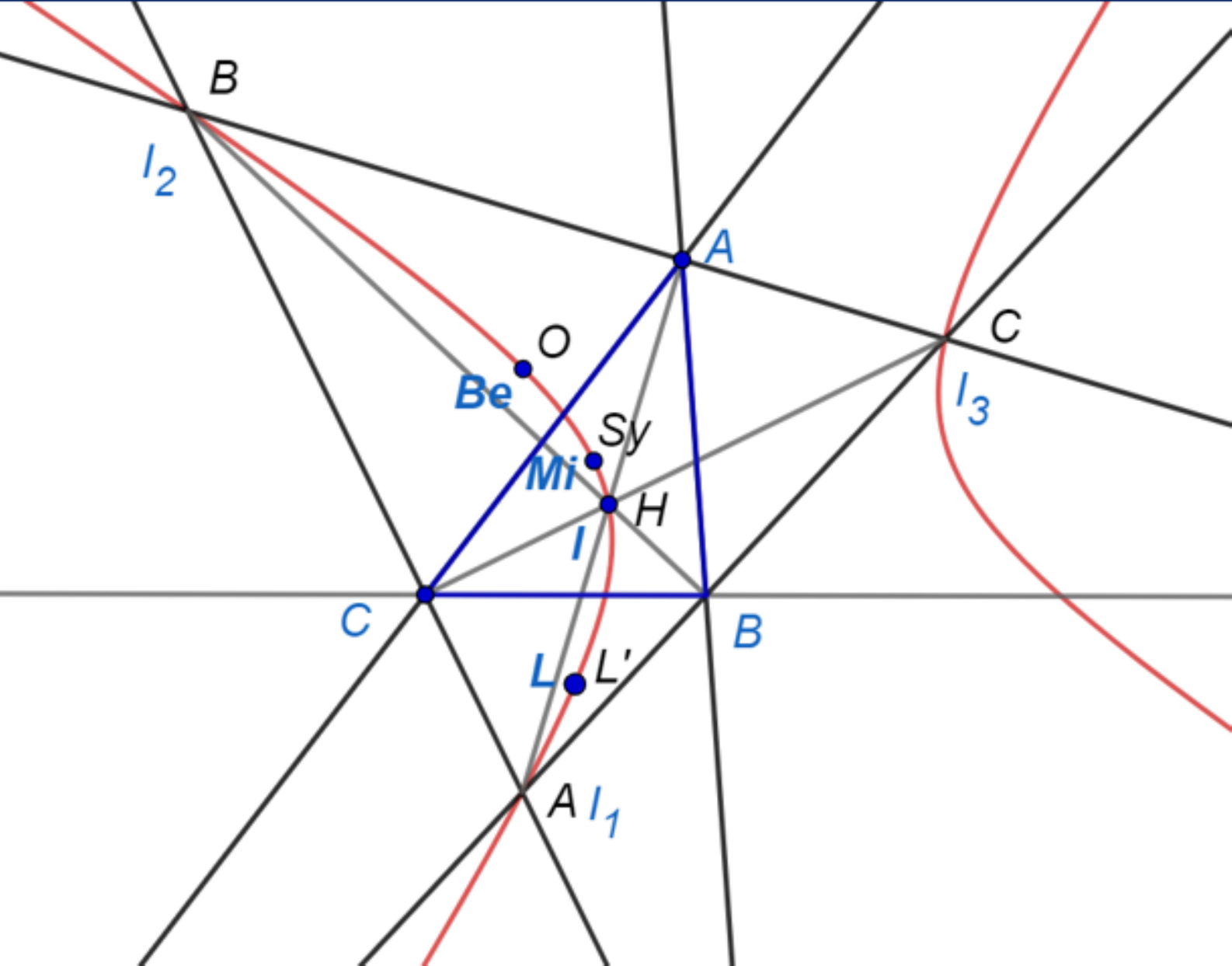}}
	\caption{New trianlge conic based on Jerabek hyperbola for the excentral triangle}
	\label{fig:JerabekNew}
\end{figure}

 \indent Therefore, we got new triangle hyperbola \ref{fig:JerabekNew}, which passes through cneters of the excircles, Bevan point, incenter, mittenpunkt and de Longchamps point. It is still rectangular as Jerabek hyperbola is. Known fact about Jerabek hypperboal is that its center is center of the Euler circle. However, Euler circle of the excentral trinangle is circumscribed circle for the base triangle. In addition, Jerabek Hyperbola is isogonaly conjugate to the Euler line. Meanwhile, Euler line for the excentral triangle is line  ($I$ (incenter), $O$ (circumcenter), $Be$ (Bevan point), $Mi'$ (isogonal conjugate for the mittenpunkt with respect to the extriangle),$ Mi''$(isogonal conjugate for the mittenpunkt with respet to the base triangle)). Therefore, we can conclude, that our new hyperbola is isogonaly conjugate to the line $(I, O, Be, Mi', Mi'')$ and its center is the circumcenter.

\begin{thm}
	New hyperbola passes though excenters, Bevan point, incenter, mittenpunkt, de Longchaps point. It is isogonal conjugate to the line  $(I, O, Be, Mi', Mi'')$, and its center is circumcircle.
\end{thm}

\indent Similarly we studied Thomson cubic for the base trianlge and matched its points with ones in the excetral trianlge.

\begin{defn}
	Thomson cubic is a curve that passes through verticies of the trianlge, middles of triangle sides, centers of the excircles, incenter, centroid, circumcenter, Lemoine point,mittenpunkt, isogonal conjugate of the mittenpunkt \cite{Klimberg}.
	\end{defn}
 \indent By apllying correspondence table between points of the base and excentral triangles \ref{tabular:table1} to the points of Thomson cubic we obtain a new triangle cubic\ref{tabular:table3}.

\begin{table}[h]
	\begin{center}
		\begin{tabular}{ | l | l |}
			\hline
			\textbf{Thomson cubic for the base triangle} & \textbf{New cubic for the excentral triangle}\\\hline
			$A, B, C$ (vertices of the base trianlge) & $H_{1}, H_{2},H_{3}$ (bases of the altitudes) \\\hline
			$M_{a}, M_{b}, M_{c}$ (middles of the base triangle sides)&$M_{ha}, M_{hb}, M_{hc}$ (midles of orthic triangle's sies)\\\hline
			$I_{1}, I_{2}, I_{3}$ (excenters) & $A, B, C$ (vertices)\\\hline
			$I$ (incenter) & $H$ (orthocenter)\\\hline
			$M$ (centroid) & $ M(H_{1}H_{2}H_{3})$ (centroid in the orthic triangle)\\\hline
			$O$ (circumcenter) & $E$ (nine-point center)\\\hline
			$Sy$ (Lemoine point) & $Sy(H_{1}H_{2}H_{3})$ (Lemoine point of the orthic triangle)\\\hline
			$Mi$ (mittenpunkt) & $Sy$ (Lemoine point)\\\hline
			\parbox[c]{6cm}{$Mi''$ (isogonal conjugate of the mittenpunkt with respet to the excentral triangle)} & \parbox[c]{6cm}{$GOT$ (gomotetic center of the orthic anf tngent triangles)}\\\hline
			
		\end{tabular}
	\end{center}
	\indent
	\label{tabular:table3}
	\caption{Mtaching points for Thomson cubic}
\end{table} 	

\indent  According to the table, we got new cubic.

\begin{figure}[h]
	\center{\includegraphics[scale=0.7]{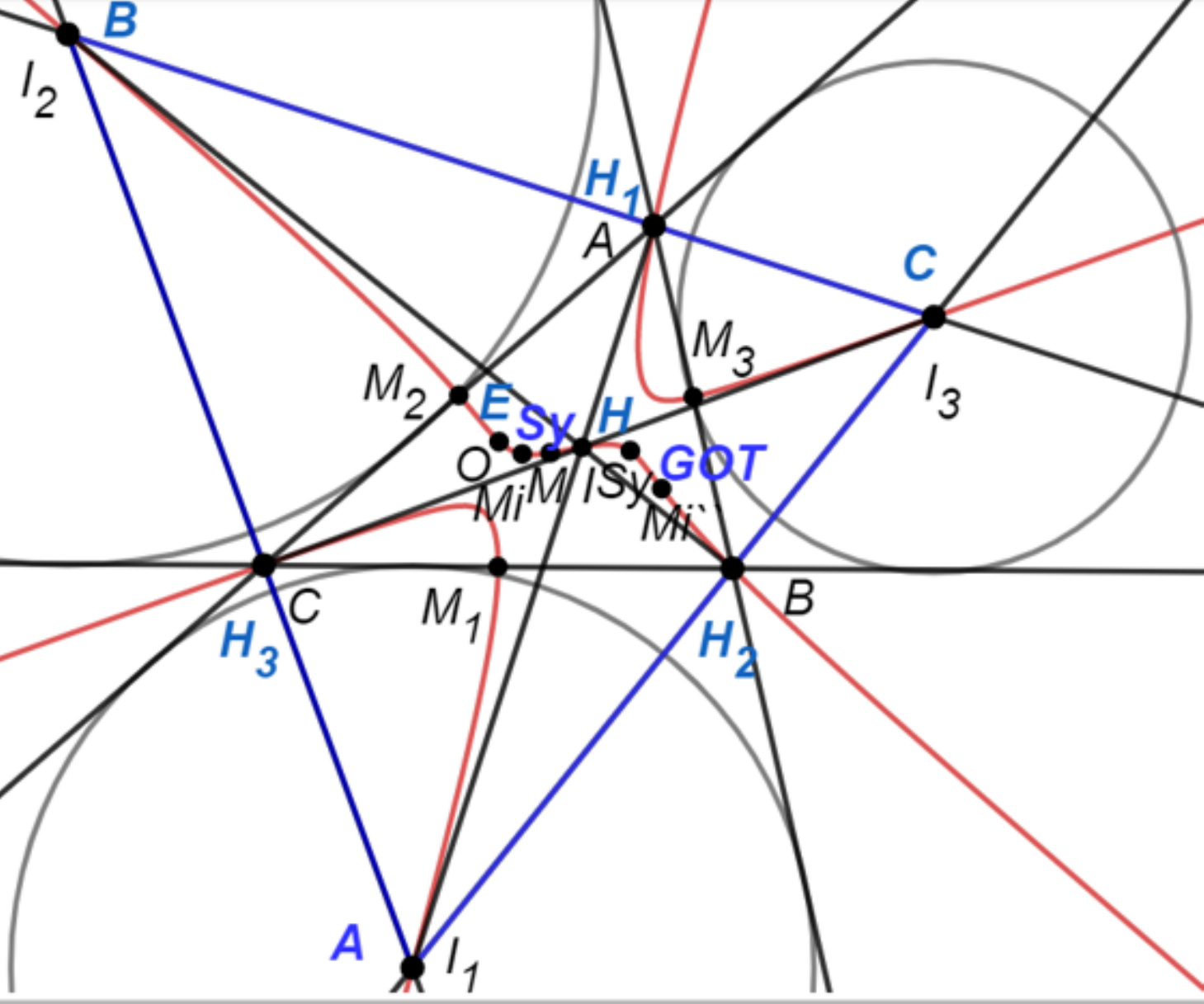}}
	\caption{New trianlge conic based on Thomson hyperbola for the excentral triangle}
	\label{fig:ThomsonNew}
\end{figure}

\begin{thm} New triangle cubic \ref{fig:ThomsonNew} passes through vertices of the triangle, bases of the altitudes, middles of the triangle sides in the orthic triangle, orthocenter, Euler point, centroid in orthic triangle, Lemoine point, and gomotetic center of the orthic and tangent triangles.
	\end{thm}

 \indent Analogically was derived new trianlge cubic based on the Darboux cubic and correspondence of its points with triangle centers in the excentral triangle.

 \begin{defn}
Darboux cubic is a curve that passes through vertices of the triangle, centers of the excircles, incenter, circumcenter, Bevan point \cite{Cundy}.
 \end{defn}

 Trianlge centers of the Darboux cubic in the base triangle were matched with points in the excentral triangle.

\begin{table}[h]
\begin{center}
\begin{tabular}{ | l | l |}
\hline
\textbf{Darboux cubic for the base triangle} & \textbf{New cubic for the excentral triangle} \\\hline
$A, B, C$ (vertices of the base triangle) & $H_{1}, H_{2}, H_{3}$ (bases of the altitudes)	\\\hline
$I_{1}, I_{2}, I_{3}$ (excenters) & $A, B, C$ (vertices) \\\hline
$I$ (incenter) & $H$ (orthocenter)\\\hline
$O$ (circumcenter) & $E$ (nine-point center)\\\hline
$Be$ (Bevan point) & $O$ (circumcenter)\\\hline
\end{tabular}
\end{center}
\indent
\label{tabular:table4}
\caption{Mtaching points for Darboux cubic}
\end{table}

\begin{figure}[h]
	\center{\includegraphics[scale=0.6]{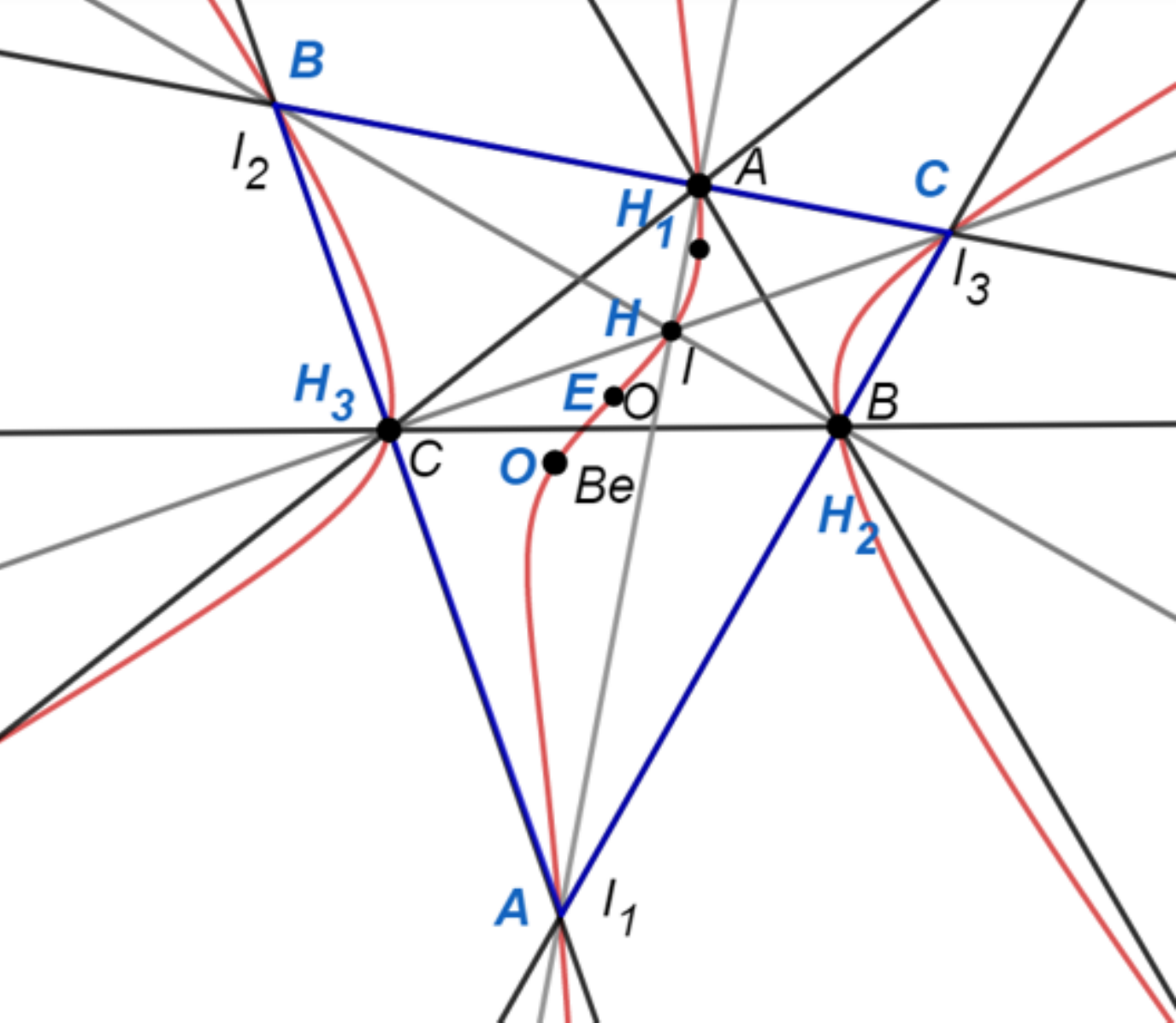}}
	\caption{New trianlge conic based on Darboux cubic for the base triangle in correspondence with excentral trianlge}
	\label{fig:DarbouxNew}
\end{figure}

 \begin{thm}
 	Therefore, we got new cubic \ref{fig:DarbouxNew}, which passes through vertices of the triangle, bases of the altitudes, orthocenter, Euler center, and circumcenter.
 \end{thm}

\indent The discussed above results were obtained from considering excentral triangle, its triangle centers and correspondence between points in the excentral and basic triangle. As a result, were derived three new triangle curves, which were not discovered before. However, to get wider results were applied the same idea to other triangles. Namely was consider medial triangle.

\begin{defn}
	Medial triangle is a triangle with vertices in the middles of the base traingle sides.
\end{defn}

\begin{table}[h]
	\begin{center}
		\begin{tabular}{ | l | l |}
			\hline
\textbf{Points in the medial triangle}& \textbf{Point in the base triangle}\\\hline
$I$ (incenter) & $Sp$ (Speaker point)\\\hline
$M$ (centroid) & $M$ (centroid)\\\hline
$O$ (circumcenter) & $E$ (nine-point center)\\\hline
$H$ (orthocenter) & $O$ (circumcenter) \\\hline
$L$ (de Longchamps point) & $H$ (orthocenter) \\\hline
$Be$ (Bevan point) &\parbox[c]{6cm}{ $Be(M_{1}M_{2}M_{3})$ (Bevan point of the medial triangle)}\\\hline
$Na$ (Nagel point) & $I$ (incenter)\\\hline
$G$ (Gergonne point) & $Mi$ (mittenpunkt)\\\hline
\parbox[c]{6cm}{$Sy_{A}$ (Lemoine point of the anticomplemetary triangle)} & $Sy$ (Lemoine point)\\\hline
$B_{3}$ (third Brocard point) & $M_{B}$ (Brocard midpoint)\\\hline
\end{tabular}
\end{center}
\indent
\label{tabular:table5}
\caption{Correspondence table between points of the medial and base triangles}
\end{table}

\indent In the same way as before, was proven the fact that some points in the medial triangle match with some points in the base triangle\ref{tabular:table5}. Proof of the mentioned facts relies on the patterns of the Euclidean geometry, some of the correspondence were proved before \cite{Klimberg}.

\begin{defn}
	Yff hyperbola is a triangle curve wich passes through centroid, orthocenter, circumcenter, and Euler center.
\end{defn}

 \indent We have considered Yff hyperbola for the base triangle and matched its point with triangle centers of the medial triangle, applying correspondence table\ref{tabular:table5}.

 \begin{table}[h]
 	\begin{center}
 		\begin{tabular}{ | l | l |}
 			\hline
 		     \textbf{Yff hyperbola for the base triangle} & \textbf{New hyperbola for the medial traingle}\\\hline
 		     $M$ (centroid) & $M$ (centroid)\\\hline
 		     $H$ (orthocenter) & $L$ (de Longchaps point)\\\hline
 		     $O$ (circumcenter) & $H$ (orthocenter)\\\hline
 		     $E$ (nine-point center) & $O$ (circumcenter)\\\hline
 		\end{tabular}
 	\end{center}
 	\indent
 	\label{tabular:table6}
 	\caption{Mtaching points for the Yff hyperbola and medial triangle}
 \end{table}

\indent We got new conic, which has vertices in the centroid and de Longchaps point, focus in the orthocenter. Directix of the Yff hyperbola is perpendicular to the Euler line and passes through center of the Euler circle. Euler line for the medial and base triangles coincide. However, center of the Euler circle of the base triangle is circumcenter for the medial. Therefore, directrix of the new hyperbola is perpendicular to the Euler line and passes through circumcenter.

\begin{figure}[h]
	\center{\includegraphics[scale=0.35]{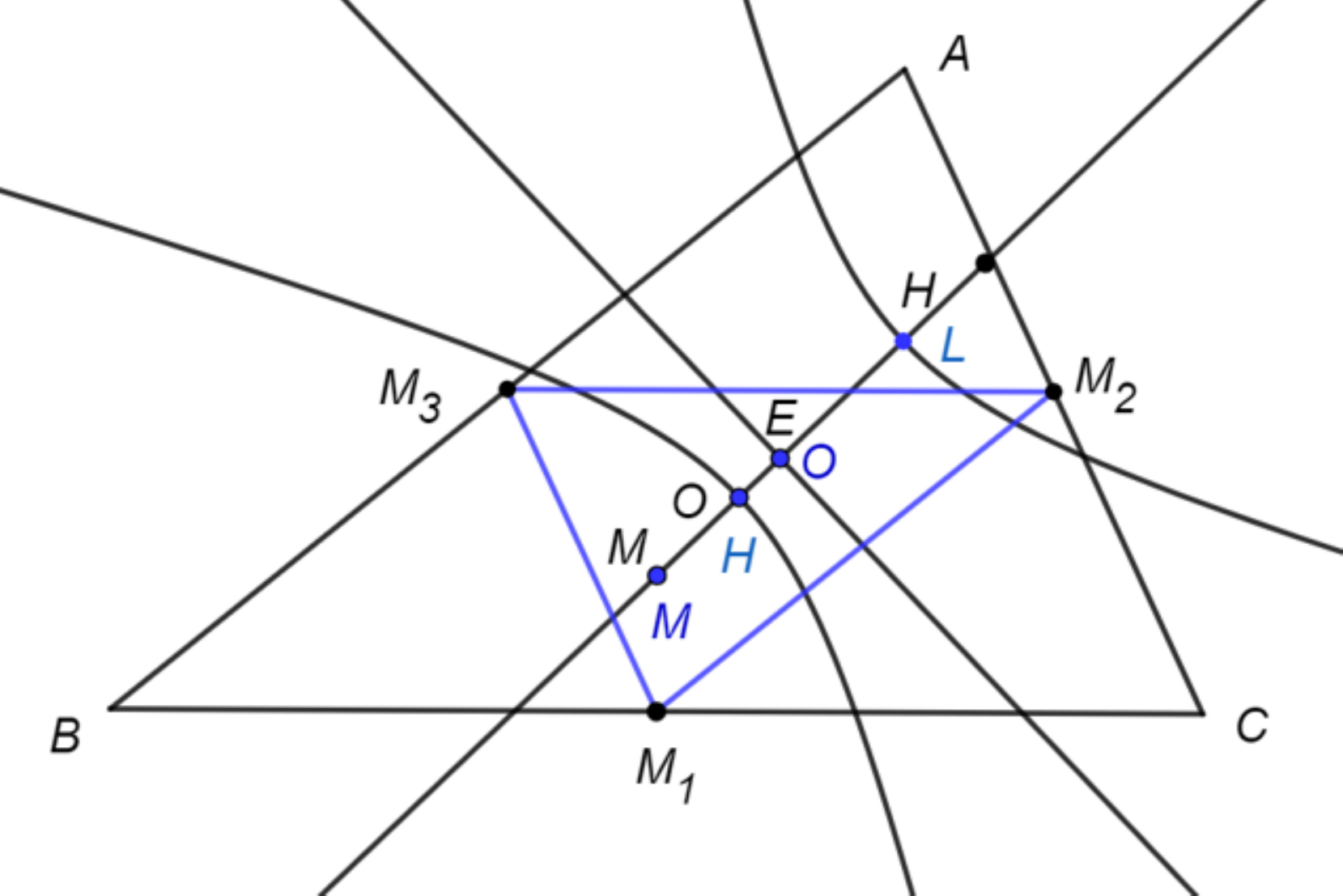}}
	\caption{New conic based on Yff hyperbola}
	\label{fig:YffNew}
\end{figure}

\begin{thm} New conic \ref{fig:YffNew} is a curve that has vertices in the centroid and de Longchaps point, focus in the orthocenter. Directrix of the new hyperbola is perpendicular to the Euler line and passes through circumcenter.
\end{thm}

\indent Let`s go further in our research and create more triangle cubics with the help of correspondence between triangle centers in medial and base triangles.

\begin{table}[h]
	\begin{center}
		\begin{tabular}{ | l | l |}
			\hline
			\textbf{Darboux cubic for the medial triangle} & \textbf{New cubic for the base triangle}\\\hline
			$A, B, C$ (triangle vertices) & $M_{1}, M_{2}, M_{3}$ (middles of the triangle sides)\\\hline
			$A_{1}, B_{2}, C_{3}$ (antipods of the triangle) & Antipods of the medial triangle\\\hline
			$I$ (incenter) & $Sp$ (Speaker point)\\\hline
			$O$ (circumcenter) & $E$ (nine-point center)\\\hline
			$H$ (orthocenter) & $O$ (circumcenter)\\\hline
			$L$ (de Longchamps point) & $H$ (orthocenter)\\\hline
			$L'$ (isogonal conjugate of the de Longchaps point) & $H_{A}$ (complementary conjugate of the orthocnter)\\\hline
		\end{tabular}
	\end{center}
	\indent
	\label{tabular:table7}
	\caption{Mtaching points for the Darboux cubic and medial triangle}
\end{table} 	

\indent We obsereved the transformation of the Darboux cubic and under the correspondence. It led us to a new cubic.

\begin{figure}[h]
	\center{\includegraphics[scale=0.6]{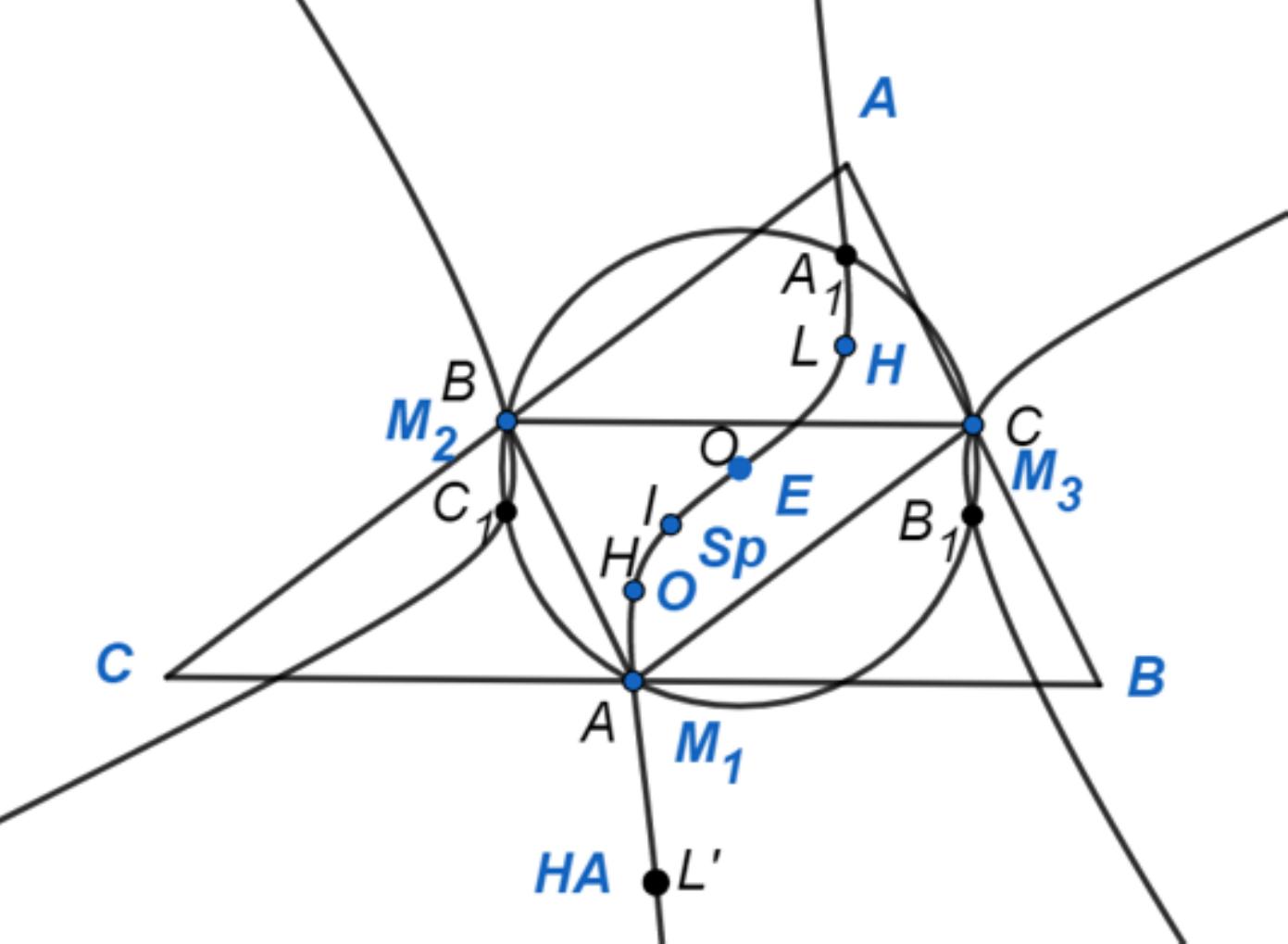}}
	\caption{New cubic based on Darboux cubic for the medial triangle}
	\label{fig:DarbouxNew2}
\end{figure}

\begin{thm}
We got new cubic \ref{fig:DarbouxNew2}, which passes through Speaker point, center of the Euler circle, circumcenter, orthocenter, complementary conjugate of the orthocenter, middles of the triangle side, and antipodes of the medial triangle.
\end{thm}

Simiraly, we take Lucas cubic for the base cubic and match it with points of the medial triangle.

\begin{defn}
 Lucas cubic is such curve which passes through triangle vertices, orthocenter, Gergone point, centroid, Nagel point, Lemoine point of the anticomplementary tringle, and vertices of th anticomplementary triangle.
\end{defn}

\indent We make the correspondence between points of the Lucas cubic in the medial triangle with triangle centers in the base triangle. As a result we obtain the following table.

 \begin{table}[h]
 	\begin{center}
 		\begin{tabular}{ | l | l |}
 			\hline
\textbf{Lucas cubic for the medial triangle} & \textbf{New cubic for the base triangle}\\\hline
$Sy_{A}$ (Lemoine point of the anticomplementary triangle) & $Si$ (Lemoine point)\\\hline
$M$ (centroid) & $M$ (centroid)\\\hline
$H$ (orthocenter) & $O$ (circumcenter) \\\hline
$Ge$ (Gergonne point) & $Mi$ (mittenpunkt) \\\hline
$Na$ (Nagel point) & $I$ (incenter) \\\hline
$L$ (de Longchamps point) & $H$ (orthocenter)\\\hline
 		\end{tabular}
 	\end{center}
 	\indent
 	\label{tabular:table8}
 	\caption{Mtaching points for the Lucas cubic and medial triangle}
 \end{table}

\begin{figure}[h]
 	\center{\includegraphics[scale=0.5]{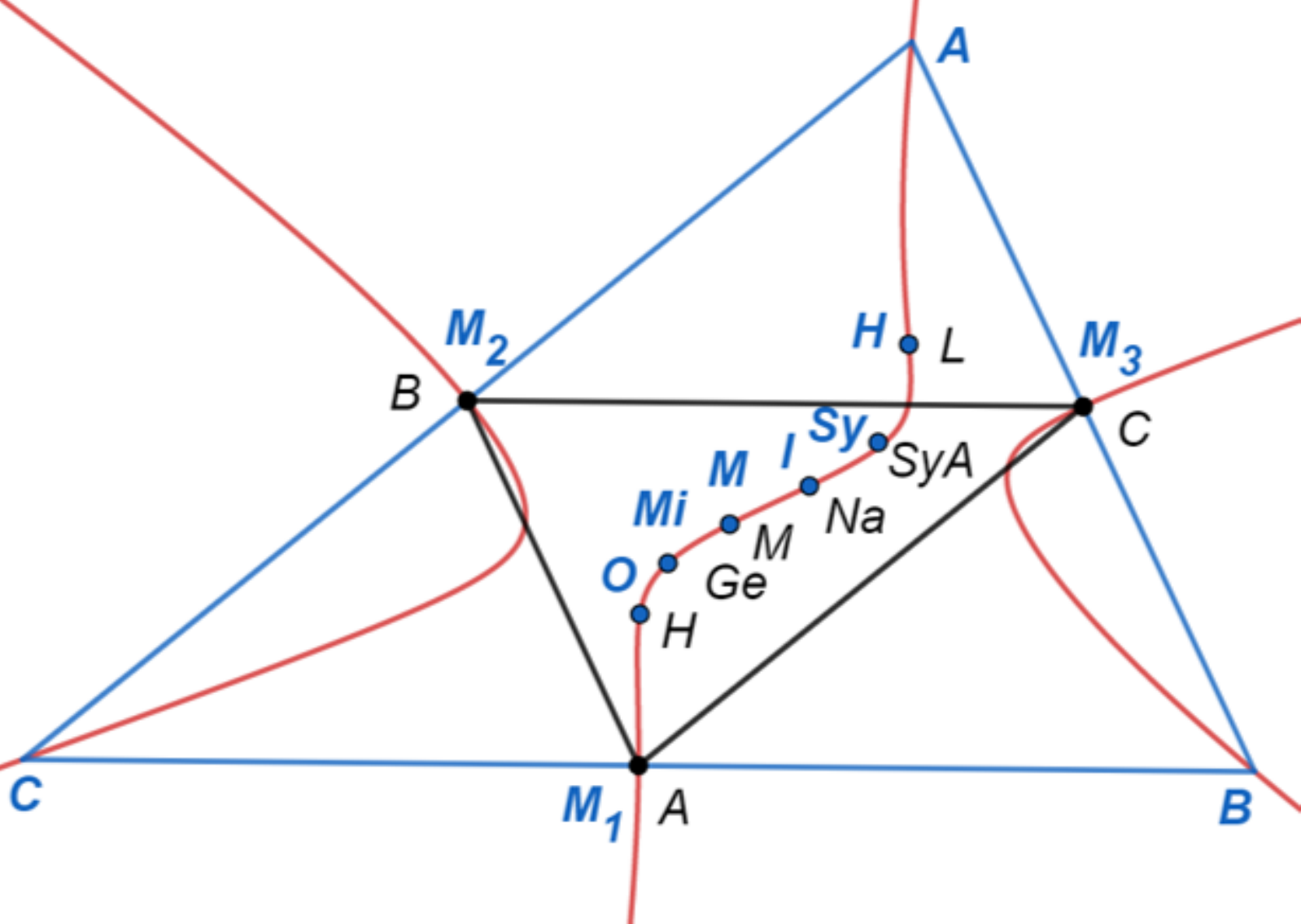}}
 	\caption{New cubic based on Lucas cubic for the medial triangle}
 	\label{fig:LucasNew}
 \end{figure}

\begin{thm}
	 New cubic \ref{fig:LucasNew}  which passes through Lemoine point, centroid, circumcenter, mittenpunkt, incenter, orthocenter, triangle vertices, and middles of the triangle sides.
\end{thm}

\indent Therefore, while applying correspondence method to the medial triangle we derived one new conic and two new cubics. In addition, we obsereved Euler and mid-arc triangles as correspondence base.

\begin{defn}
	Euler triangle is triangle with vertices in the intersection points of the triangle altitudes and nine-point circle.
	\end{defn}
\begin{defn}
	Mid-arc triangle is a triangle with vertices in the middles of the arcs of the circumcircle.
	\end{defn}

\indent Let's firstly consider correspondence of points between Euler and base triangles.

\begin{table}[h]
	\begin{center}
		\begin{tabular}{ | l | l |}
			\hline
		\textbf{Points in the Euler triangle} & \textbf{Points in the base triangle}\\\hline
		$I$ (incenter) & $M_{IH}$ (midpoint of incener and orthocenter)\\\hline
		$M$ (centroid) & $M_{MH}$ (midpoint of centroid and orthocenter)\\\hline
		$O$ (circumcenter) & $E$ (nine-point center)\\\hline
		$H$ (orthocenter) & $H$ (orthocenter) \\\hline
		$N$ (Nagel point) & $F$ (Furhman point) \\\hline
		$L$ (de Longchamps point) & O (circumcenter)\\\hline
		\end{tabular}
	\end{center}
	\indent
	\label{tabular:table9}
	\caption{Correspondence table between points in the Euler ad base triangles}
\end{table}

\indent  According to the above table we have buil the correspondence between points of the Darboux cubic for the Euler triangle and triangle centers of the base triange.

\begin{table}[h]
	\begin{center}
		\begin{tabular}{ | l | l |}
	\hline
     \textbf{Darboux cubic in the Euler triangle} & \textbf{New cubic for the base triangle}\\\hline
     $A, B, C$ (vertices) & $E_{1}, E_{2}, E_{3}$ (vertices of the Euler triangle)\\\hline
     $A', B', C'$ (antipods of the triangle vertices) & $M_{1}, M_{2}, M_{3}$ (middles of the triangle sides)\\\hline
     $I$ (incenter) & $M_{IH}$ (midpoint of incenter and orthocenter)\\\hline
     $O$(circumcenter) & $E$ (nine-point center)\\\hline
     $H$ (orthocenter) & $H$ (orthocenter)\\\hline
     $L$ (de Longchamps point) & $O$ (circumcenter)\\\hline
		\end{tabular}
	\end{center}
	\indent
	\label{tabular:table10}
	\caption{Mtaching points for the Darboux cubic and Euler triangle}
\end{table}

\begin{figure}[h]
	\center{\includegraphics[scale=0.4]{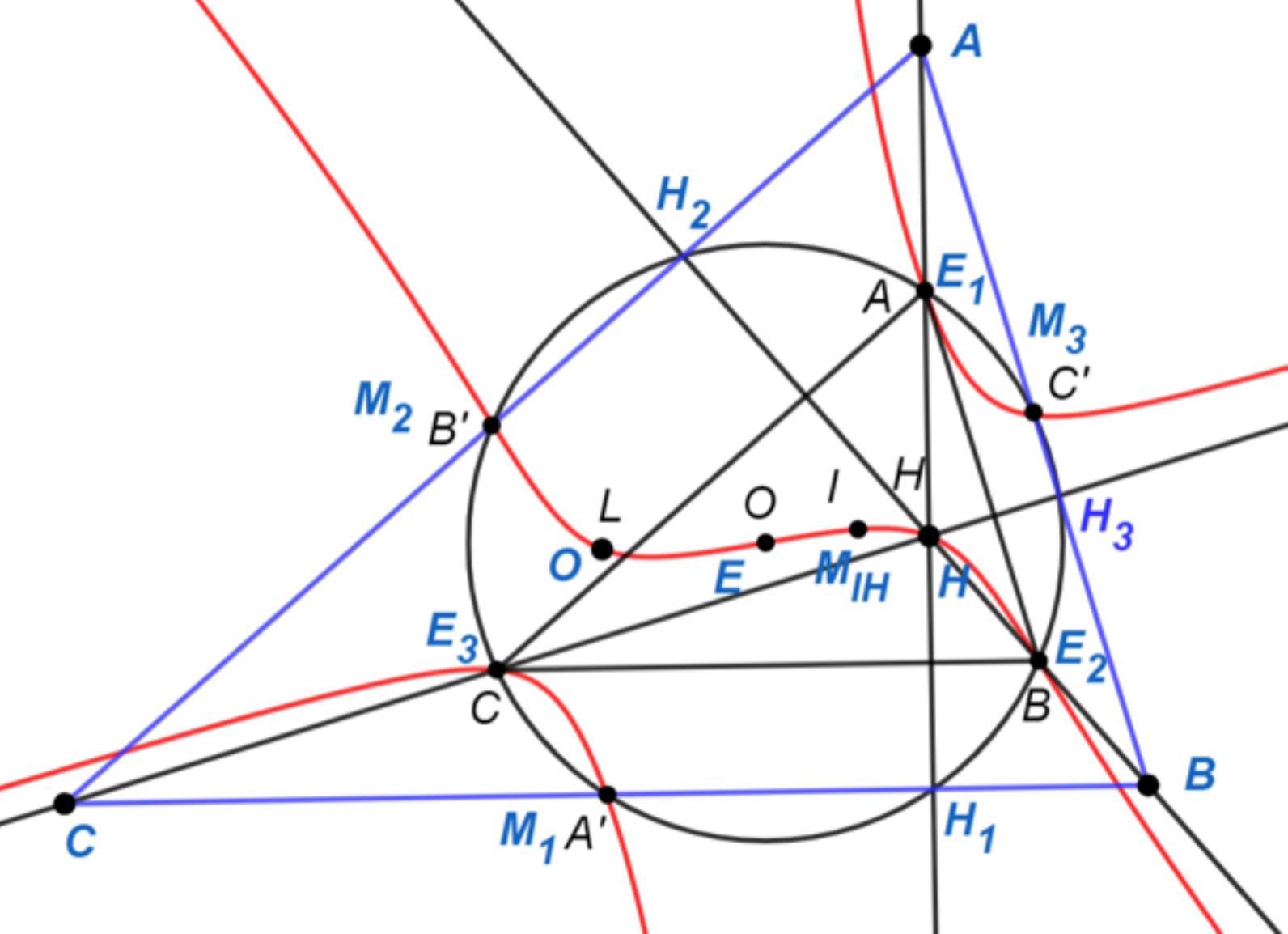}}
	\caption{New cubic based on Darboux cubic in the Euler trianle}
	\label{fig:DarbouxNew3}
\end{figure}

\begin{thm}
	 New cubic \ref{fig:DarbouxNew3} passes through circumcenter, orthocenter, Euler center, midpoint of the incenter and the orthocenter, vertices of the Euler triangle and middles of the triangle sides.
\end{thm}

\indent Finally, we observe correspondence of traingle centers between mid-arc and base trinalge.

\begin{table}[h]
	\begin{center}
		\begin{tabular}{ | l | l |}
			\hline
			\textbf{Points for the mid-arc triangle} & \textbf{Points for the base triangle}\\\hline
		$O$ (circumcenter) & $O$ (circumcenter)\\\hline
		$H$ (orthocenter) & $I$ (incenter) \\\hline
		$Sy$ (Lemoine center) &\parbox[c]{6cm}{ $M_{MiI}$ (midpoint of mittenpunkt and incenter)}\\\hline
		$L$ (de Longchamps point) & $Be$ (Bevan point)\\\hline
		$K$ (Kosnita point) & $S$ (Schiffler point)\\\hline
		\end{tabular}
	\end{center}
	\indent
	\label{tabular:table11}
	\caption{Correspondence table between points in the mid-arc and base triangles}
\end{table}

\indent Based on the correspondence between point sof the mid-arc and base triangles we discovered new cubic which is based on Jerabek hyperbola.

 \begin{table}[h]
 	\begin{center}
 		\begin{tabular}{ | l | l |}
 			\hline
 		\textbf{Jerabek hyperbola for mid-arc triangle} & \textbf{New hyperbola for the base triangle}\\\hline
 		$A, B, C$ (vertices)& \parbox[c]{6cm}{$A_{1}, A_{2}, A_{3}$ (middles of the arcs of the circumcircle)}\\\hline
 	    $O$ (circumcenter) & $O$ (circumcenter)\\\hline
 	    $H$ (orthocenter) & $I$ (incenter) \\\hline
 	    $Sy$ (Lemoine point) & \parbox[c]{6cm}{$M_{MiI}$ (midpoint of mittenpunkt and inenter)}\\\hline
 	    $K$ (Kosnita point) & $S$ (Schiffler point) \\\hline
 	    \parbox[c]{6cm}{$L'$ (isogonal conjugate of the de Longchamps point)} & \parbox[c]{6cm}{$Be'$ (isogonal conjugate of the Bevan point)}\\\hline
 		\end{tabular}
 	\end{center}
 	\indent
 	\label{tabular:table12}
 	\caption{Mtaching points for the Jerabek hyperbola and Euler triangle}
 \end{table}

\begin{figure}[h]
 	\center{\includegraphics[scale=0.5]{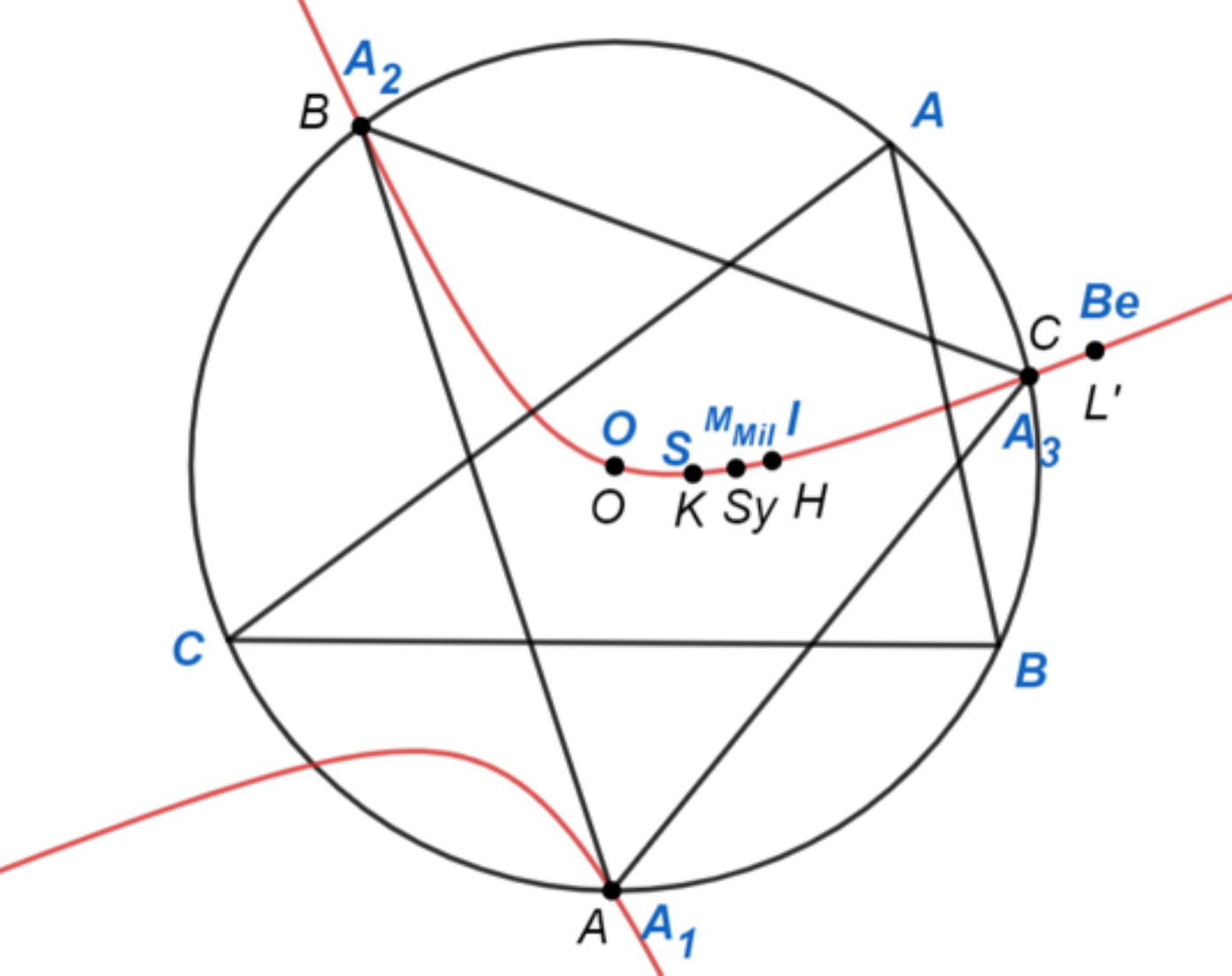}}
 	\caption{New cubic based on the Jerabek hyperbola for the mid-arc triangle}
 	\label{fig:JerabekNew2}
 \end{figure}

 \begin{thm}
 	New conic \ref{fig:JerabekNew2} is rectangular and passes through circumcenter, incenter, midpoint of mittenpunkt and incenter, Schiffler point, and isogonaly conjugate point to the Bevan point.
 	\end{thm}

\indent Therefore, during the research of the triangle curves were derived three new triangle conics nad five new triangle cubics. This is a significant result and leaves room for new investigations.

\indent Since, geometry of conic sections and other triangle curves are broadly used in the projective geometry we looked on the obtained result through the prism of the projective geometry.

\indent According to the Pascal's theorem if six arbitrary points are chosen on a conic and joined by line segments in any order to form a hexagon, then the three pairs of opposite sides of the hexagon meet at three points which lie on a straight line.

\indent Let`s consider the first derived triangle curve based on Jerabek hyperbola for the excentral triangle. New hyperbola passes though excenters, Bevan point, incenter, mittenpunkt, de Longchaps point. We bulit a hexagon with verices in the given triangle centers and apply Pascal's theorem.

\indent Let $I_{1}, I_{2}$ be excenters, and $Be$, $Mi$, $L$ be Bevan point, $Mi$ mittenpunkt, de Longchamps point, respectively. We get the following results:
 \begin{cor}
  Concurent points of $I_{2}Be$ and $LI$, $BeMi$ and $I_{1}L$, $MiI$ and $I_{2}I_{1}$ belong to one line.
 \end{cor}
 \begin{cor}
 	 Concurent points of segments $I_{2}Mi$ and $BeI_{1}$, $BeL$ and $II_{2}$, $MiL$ and $II_{2}$ lie on one line.
 \end{cor}

\indent Simiraly, we have applied the same idea for the hexagon inscribed in the new hyperbola derived from the Jerabek hyperbola for the mid-arc triangle.

 \indent Let $A_{2}, A_{3}$ be middles of the arcs of the circumcircle, and $Be$, $I$,  $S$,  $O$ be Bevan point, incenter, Speaker point, circumcenter, respectively. The following facts were discovered:

\begin{cor}
	Points of intersection of lines $A_{2}Be$ and $SI$, $IA_{3}$ and $OA_{2}$, $BeA_{3}$ and $OS$ belong to one line.
\end{cor}

\begin{cor}
	Points of intersection of line segments $BeO$ and $A_{3}A_{2}$, $BeS$ and $IA_{2}$, $A_{3}S$ and $IA_2$ belong to a straight line.
\end{cor}

Moreover, combination of two of the discovered triangle cubics gives us very interesting corollary as well. Let's consider new cubic derived from the Darboux cubic for the excentral triangle and new cubic constructed with the base Darboux cubic with respect to the medial triangle. The first mentioned new cubic passes through bases of altitudes, vertices, orthocenter, nine-point center, circumcenter, let's name it $P(x,y)$. The second mentioned new cubic passes through middles of the triangle sides, Speaker point, nine-point center, circumcenter, orthocenter, let`s name it $Q(x,y)$. We may notice that this two cubics pass through three common points which are nine-point center, circumcenter and orthocenter. Moreover, this three points belong to Euler line, let it has an equation $ax+by+c=0$. Since, we have two cubics which pass through points which belong to one line, there exists such integer $t$ such that the following holds:
$P(x,y)-tQ(x,y)$ \vdots\ $ax+by+c$. Therefore, Euler line is linear component of the composition of two new cubics. In addition, points of intersection of the linear component with the curve are inflection points \cite{walker}.

\begin{cor}
	Euler line is a linear component of the composition of new triangle cubic (passes through bases of altitudes, vertices, orthocenter, nine-point center, circumcenter) and new triangle cubic (passes through middles of the triangle sides, Speaker point, nine-point center, circumcenter, orthocenter). Moreover, orthocenter, circumcenter, and nine-point center are inflection point of the composition of these two curves.
\end{cor}

\indent The above corollaries prove that the discovered in the research new triangle curves could be applied in different geometric areas and studied in advanced.

\begin{rem} A further continue of our research consists in the same analysis of singularities as provided by second author in \cite{SkCubic, SkOs} for cubic obtained by us in the presented work.
\end{rem}
\indent

\indent
\section{Conclusion}

 
 During the research were discovered three new triangle conics and five new triangle cubics, what is very significant result for the classical geometry. In addition, was shown that proceedings of the study could be applied not only in Euclidian space, but in projective as well. However, the main result of the research was developed algorithm of deriving new triangle curves. This opens an opportunity for creating more triangle curves, while applying the method for various triangles, points, and geometric constructions.

\indent The developed idea significantly simplifies the question of creating curves passing through triangular centers. However, it opens up a number of new questions. Which interesting properties do new curves have? What is the topological nature of these transformations? Is it possible to apply a similar idea to non-Euclidean objects? Could one use the same method over an arbitrary finite field? Can this idea be further generalized?

\end{document}